\documentclass[final]{jcmlatex}
\usepackage{array, boldline, makecell, booktabs}
\usepackage{cases}
\usepackage{epstopdf}
\usepackage{arydshln}
\usepackage{xcolor}
\usepackage{graphics}
\usepackage{graphicx}
\usepackage{amssymb}
\usepackage{algorithm}
\usepackage{algorithmic}
\usepackage{longtable}
\usepackage[colorlinks,bookmarksopen,bookmarksnumbered,citecolor=blue,urlcolor=red]{hyperref}
\usepackage[figuresright]{rotating}

\begin{document}

\def\DOI{?}
\def\JCMvol{}
\def\JCMno{}
\def\JCMyear{}
\def\JCMreceived{?}
\def\JCMrevised{?}
\def\JCMaccepted{?}
\def\argmin{\mbox{argmin}}
\setcounter{page}{1}
\newtheorem{lem}{Lemma}[section]
\newtheorem{dingyi}{Definition}[section]
\newtheorem{exa}{Example}[section]
\newtheorem{coro}{Corollary}[section]
\newtheorem{rem}{Remark}[section]
\newtheorem{thm}{\bf Theorem}
\newcommand{\aaa}{\left \langle}
\newcommand{\bbb}{\right \rangle}
\newcommand{\ccc}{\left (}
\newcommand{\ddd}{\right )}
\newcommand{\eee}{\Big \{}
\newcommand{\fff}{\Big \}}
\newcommand{\C}{\mbox{$\mathbb{C}$}}
\newcommand{\R}{\mbox{$\mathbb{R}$}}
\newcommand{\SP}{\mbox{$\mathcal{SP}$}}
\newcommand{\RR}{\mbox{$\mathcal{R}$}}
\newcommand{\PP}{\mbox{$\mathbf{P}$}}
\newcommand{\T}{\mbox{$\mathbf{T}$}}
\newcommand{\0}{\mbox{${\bf 0}$}}
\newenvironment{breakablealgorithm}
  {
   \begin{center}
     \refstepcounter{algorithm}
     \hrule height.8pt depth0pt \kern2pt
     \renewcommand{\caption}[2][\relax]{
       {\raggedright\textbf{Algorithm~\thealgorithm} ##2\par}%
       \ifx\relax##1\relax 
         \addcontentsline{loa}{algorithm}{\protect\numberline{\thealgorithm}##2}%
       \else 
         \addcontentsline{loa}{algorithm}{\protect\numberline{\thealgorithm}##1}%
       \fi
       \kern2pt\hrule\kern2pt
     }
  }{
     \kern2pt\hrule\relax
   \end{center}
  }
\def\argmin{\mbox{argmin}}
\def\QEDclosed{\mbox{\rule[0pt]{1.3ex}{1.3ex}}} 
\def\QEDopen{{\setlength{\fboxsep}{0pt}\setlength{\fboxrule}{0.2pt}\fbox{\rule[0pt]{0pt}{1.3ex}\rule[0pt]{1.3ex}{0pt}}}} 
\def\QED{\QEDclosed} 
\def\proof{\noindent{\bf Proof}: } 
\def\endproof{\hspace*{\fill}~\QED\par\endtrivlist\unskip}

\markboth{Haochen Jiang, Dongdong Liu, Xianping Wu and Xu Yang}
         {Accelerated Kaczmarz methods via randomized sketch techniques for solving consistent linear systems}
\title{\Large Accelerated Kaczmarz methods via randomized sketch techniques for solving consistent linear systems}
\author{Haochen Jiang$^{a}$\ \ \ \ Dongdong Liu$^{a}$\footnote{\texttt{corresponding author(ddliuresearch@163.com, ddliu@gdut.edu.cn)}}\ \ \ \ Xianping Wu$^{a}$ \ \ \ \ Xu Yang$^{a}$
\affiliation{$^a$ School of Mathematics and Statistics, Guangdong University of Technology, Guangzhou, China}}

 \maketitle

\begin{abstract}
Motivated by the randomized sketch to solve a variety of problems in scientific computation, {we improve both the maximal weighted residual Kaczmarz method and the randomized block average Kaczmarz method using two new randomized sketch techniques}. Besides, convergence analyses of the proposed methods are provided. {Furthermore, we establish an upper bound for the discrepancy between the numerical solutions obtained via the proposed methods and those derived from the original approaches.} Numerical experiments demonstrate that the new methods perform better than the existing ones in terms of the running time with the same accuracy.
\end{abstract}

\section{Introduction}\label{section1}
Let $\mathbb{R}$, $\mathbb{R}^{r}$ and $\mathbb{R}^{m\times n}$ denote the real field, the set of all $r$-dimensional  column vectors and the set of all $m\times n$ matrices, respectively. Consider the following least squares problems for consistent linear systems:
	\begin{equation}\label{eq1}
	\min_{x\in {\mathbb{R}}^n}	\|A x - b\|_2,
	\end{equation}
	where $A \in \mathbb{R}^{m\times n},~b \in \mathbb{R}^{m}$ and $m \gg n$. In this paper, we will investigate the improved Kaczmarz methods with the randomized matrix to solve the linear systems \eqref{eq1}.
	
	\par
The Kaczmarz method, originally proposed by Stefan Kaczmarz \cite{kac}, is a cyclic projection algorithm for solving least squares problems \eqref{eq1}. This versatile method has widespread applications in diverse fields including image reconstruction \cite{popa,Li,kam,ram,huliu,popa1}, signal processing \cite{lor,abd,xu,wu,byr}, and distributed computing \cite{elb,pas,aved}. The Kaczmarz method is designed to address large-scale problems since only one row of the coefficient matrix is involved in each iteration. Applying the randomization technology, Strohmer and {Vershynin} \cite{str} introduced the randomized Kaczmarz method (RK) that no longer follows a certain determined order but randomly selects the updated rows according to probability. However, if the Euclidean norm of each row of the residual matrix is the same, the RK method degenerates into the standard Kaczmarz method. In order to overcome this curse, Bai and Wu \cite{bai} proposed the greedy randomized Kaczmarz method (GRK), which captures the large components of the residual vector during each iteration. Based on the GRK method, various improved methods emerged, e.g., see \cite{jj,baiz,guo,zheng}. By combining the relaxation technique, Bai and Wu \cite{bai1} gave the relaxed greedy randomized Kaczmarz method (RGRK) and the convergence factor is noted to be the smallest if the relaxation parameter  $\theta = 1$, which is called the maximal weighted residual Kaczmarz method (MWRK) by Du and Gao \cite{du}. Elfving \cite{elf} proposed a block Kaczmarz method (BK), which chooses a block uniformly at random at each iteration enforcing multiple constraints synchronously, and then Needell and Tropp \cite{nee} gave a randomized block Kaczmarz method (RBK) for solving the {overdetermined} least squares problems. The RBK method divides the rows of the coefficient matrix into several subblocks and selects the blocks using the principle of uniform sampling in the projection process. Liu and Gu \cite{gu} combined GRK method and RBK method, and then got a greedy randomized block Kaczmarz method (GRBK), which grabs the bulk of residual mode length in the iteration process. Combining the RBK method with the randomized extended Kaczmarz method (REK) \cite{zou}, Needell et al \cite{need} constructed the randomized double block Kaczmarz method (RDBK) for solving inconsistent systems. Chen and Huang\cite{chen} gave the error estimate of RDBK method. Later, utilizing $k$-means clustering with the RBK method, Jiang et al. \cite{jiang}  got a new randomized block Kaczmarz method which exploits the greedy technique by capturing large residuals of different block linear systems. Necoara \cite{nec} further gave a randomized average block Kaczmarz method (RaBK), which was applied to solve a consistent tensor system under the $t$-product structure\cite{liao}.
	\par
	Recently, by employing a count sketch matrix $\mathbf{S}$, Zhang and Li \cite{zhang} transformed the optimization problem \eqref{eq1} into the following one:
		\begin{equation}\label{eq22}
		\min_{x\in {\mathbb{R}}^n}	\|{\mathbf{S}}A x - {\mathbf{S}}b\|_2,
	\end{equation}
and proposed a count sketch maximal weighted residual Kaczmarz method. The idea of this method is to map the rows of the matrix $A$ to a lower dimension space through a random mapping and binary matrix, and then randomly flip the result of the mapping through a random diagonal matrix. This process can effectively reduce the size of the matrix involved in the calculation while preserving most information  from the original problem. However, how to select a more effective random matrix is a challenging work. In this paper, we further investigate  efficient random matrices to improve the block randomized Kaczmarz method. In the following, we list the main contributions of this paper as follows:
	\begin{enumerate}
		\item By employing two randomized sketch matrices, we propose an improved maximal weighted residual Kaczmarz method and a randomized block average Kaczmarz method. Besides, convergence analyses of the proposed methods are provided.
		\item {We give the upper bound on the deviation between numerical solutions from the proposed methods and those obtained from the original approaches.}
		\item We present some numerical experiments showing that the proposed methods perform well in terms of the running time for large scale numerical examples.
	\end{enumerate}
	\par
	The organization of this article is as follows. In Section 2, first, we introduce some notations. Second, we review the maximal weighted residual Kaczmarz method and the randomized average block Kaczmarz method. Finally, we list some existing conclusions which will be used in this paper. In Section 3, we propose new randomized Kaczmarz methods employing the new random sketching matrices. The convergence of the proposed methods is given in Section 4. Numerical experiments are presented in Section 5 to illustrate the effectiveness of the new methods. The final section is the conclusion.

\section{Preliminaries}\label{section2}
    In this paper, for any random variable $\xi$, let $E[\xi]$ denote its expectation. For any vector $q \in \mathbb{R}^{m}$, $ q_i $, $ q^\ast $ and $ \| q \|_2 $ represent the $ i $-th entry, the transpose and the Euclidean norm of $ q $, respectively. For any matrix $ A \in \mathbb{R}^{m\times n} $, $ A^{(i)} $, $ A^\ast$, $ \sigma_i(A) $, $ \sigma_r(A) $, $ \| A \| _F $ represent the $i$-th row, the transpose, the $i$-th singular value, the smallest nonzero singular value and the Frobenius norm of $ A $, respectively. Besides, $ R (A) $ and $ R (A)^\bot $ represent the column space of $ A $ and its orthogonal complement, respectively. For a positive integer $m$, let the set $[m]\equiv\{1,2,\cdots,m\}$.
	\par
Next  we recap the maximal weighted residual Kaczmarz method \cite{du} given in Algorithm \ref{alg:Framwork1}.
	\begin{algorithm}[htb]
		\caption{ \cite{du} The maximal weighted residual Kaczmarz method(MWRK)}
		\label{alg:Framwork1}
		\begin{algorithmic}[1]
			\REQUIRE ~~ $A,b~\text{and}~x_0$\\	
			\STATE For $k=0, 1, 2, \cdots $ do until satisfy the stopping criteria
			\STATE Compute $i_k = \arg \max\limits_{1\leq i_k\leq m}\{\frac{\mid b^{(i_k)}-A^{(i_k)}x_k\mid^2}{\parallel A^{(i_k)}\parallel _2 ^2} \}$;
			\STATE Set $x_{k+1} = x_k+\frac{b^{(i_k)}-A^{(i_k)}x_k}{\parallel A^{(i_k)}\parallel_2^2}(A^{(i_k)})^\ast$;
			\STATE end for
				\RETURN $x_{k+1}$;
		\end{algorithmic}
	\end{algorithm}
Compared with the RGRK method, the row index selection strategy used in the MWRK method is simpler. In addition, the MWRK method operates without requiring an index set. Due to less computational cost of the row index selection strategy, the MWRK method is more efficient than the RGRK method. Furthermore, there are also some further improvements, e.g., see the works of \cite{zhang} and \cite{zhang1}. The randomized average block Kaczmarz method \cite{nec} is given in Algorithm \ref{alg:Framwork2}.
	\begin{algorithm}[htb]
		\caption{ \cite{nec} The randomized average block Kaczmarz method(RaBK)}
		\label{alg:Framwork2}
		\begin{algorithmic}[1]
			\REQUIRE ~~	$ A,b,x_0,\{\alpha_k\}_{k\ge 0}~\text{and}~\{\omega_k\}_{k\ge0} $\\		
		    \STATE For $k=0, 1, 2, \cdots $ do until satisfy the stopping criteria
			\STATE Draw sample $J_k\sim \textup{P}$ and update:
			\STATE  $x_{k+1} = x_k-\alpha_k(\sum\limits_{i\in J_k}\omega_k^i\frac{A^{(i)} x_k-b^{(i)}}{\|A^{(i)}\|_2^2}(A^{(i)})^\ast)$;
			\STATE end for
				\RETURN $x_{k+1}$;
		\end{algorithmic}
	\end{algorithm}	
Here $J_k=\{i_k^1,\cdots, i_k^{\tau_k}\}\subseteq[m]$ is the set of indexes corresponding to the rows selected at the $k$-th iteration of the size $\tau_k\in [1,m]$, {and} $\textup{P}$ denotes the probability distribution over the collection of subsets of indexes of $[m]$. Assume that the probability $\textup{P}$ is chosen in such a way that for any index $i\in [m]$, the probability of $i$ {in $J$ satisfies $p_i>0$}. Besides, $\{\omega_k\}_{k=0}^{\infty}$ represents the weights sequence which satisfies $\omega_k^i\in[0,1]$ and $\sum\limits_{i\in J_k}\omega_k^i=1$. For an index set $J\subset[m]$, by $A_J\in\mathbb{R}^{|J|\times n}$ we denote the matrix with the rows $A_i$ for $i\in J$. The projection of a point $x\in \mathbb{R}^{n}$ onto a closed convex set $\chi\in \mathbb{R}^{n}$ is denoted by $\Pi_\chi(x)=\mathop{\arg\min}_z\{\|x-z\|\mid z\in \chi\}$. We use $\lambda_{\textup{min}}^{nz}(\cdot)$ to denote the minimum nonzero eigenvalue of a given matrix.

\begin{rem}\label{alpha}
The weights $\omega_k^i$ in Algorithm \ref{alg:Framwork2} are chosen to satisfy $0<\omega_{\textup{min}}\leq\omega_k^i\leq\omega_{\textup{max}}<1$. Let $\lambda_{\textup{max}}^{\textup{block}}=\max\limits_{J\sim \textup{P}}\lambda_{\textup{max}}(A_J^\ast \textup{diag}(\frac{1}{\|A^{(i)}\|_2^2}, i\in J) A_J)$ and $\bar{\omega}_i^k=\frac{\omega_i^k}{\|A^{(i)}\|_2^2}$, where $\textup{diag}(c_i,i\in J)$ is a $|J|\times |J|$ diagonal matrix with its $i$ diagonal entries being $c_i$ and $|J|$ denotes the  number of elements in the set $J$.

In this paper, we consider two types for determining the sequence $\{\alpha_k\}$:
\begin{description}
	\item[(i)] The constant stepsize $\alpha_k=\alpha=\frac{(2-\eta)\omega_{\textup{min}}}{\omega_{\textup{max}}^2\lambda_{\textup{max}}^{\textup{block}}}$ for all $k$.
	\item[(ii)]The adaptive stepsize $\alpha_k=(2-\eta)L_k$ where $\eta\in(0,1]$ and
	\[
	L_k=\begin{cases}
		\frac{\Sigma_{i\in J_k}\bar{\omega}_i^k(A^{(i)}x_k-b^{(i)})^2}{\|\Sigma_{i\in J_k}\bar{\omega}_i^k(A^{(i)}x_k-b^{(i)})(A^{(i)})^\ast\|_2^2} & if~\exists i\in J_k~s.t~ A^{(i)}x_k-b^{(i)}\neq 0, \\
		\frac{1}{\lambda_{\textup{max}}(A_{J_k}^\ast \textup{diag}(\bar{\omega}_i^k,i\in J_k)A_{J_k})} & otherwise.
	\end{cases}
	\]
\end{description}
\end{rem}
\par
The iteration is equivalent to a convex combination of the updates of the random Kaczmarz method with a certain step size in different directions. As for the randomized block method, it is necessary to solve the Pseudo-inverse of a sub-matrix, which is relatively time-consuming. Therefore, in order to avoid this drawback, the technique of average block was introduced and the RaBK method was further proposed \cite{nec}. 

\par
Next, we give some lemmas, which will be used for the convergence analysis of the proposed methods.

	\begin{lem} \textup{(}{\cite{meng}},~Theorem 1\textup{)}\label{lem}
		For any matrix $A\in \mathbb{R}^{m\times n}$, $\mathbf{S}$ is an $d\times m$ matrix that embeds the column space of $A$ into \(\mathbb{R}^d\) while approximately preserving the $\ell_2$ norms of all vectors in that subspace. Then we have that
		\begin{equation}\label{eq23}
			(1-\epsilon)\|Ax\|_2 \leq \|\mathbf{S}Ax\|_2 \leq (1+\epsilon)\|Ax\|_2, \text{for all}~x \in \mathbb{R}^n.
		\end{equation}
	\end{lem}


	\begin{lem}\textup{(}{\cite{wood},~Theorem 2.6}\textup{)}\label{lem1}
		 Let $\delta,\epsilon\in{(0,1)}$ and $\mathbf{S}\in\mathbb{R}^{d\times m}$ be a sparse embedding matrix with $d=O(n^2/(\delta\epsilon^2))$. Then for a given matrix $A\in \mathbb{R}^{m \times n}$, $\mathbf{S}$ is a $(1\pm\epsilon)$ $\ell_2$-subspace embedding of for $A$ with probability $1-\delta$. 
	\end{lem}

	\begin{lem} \textup{(}{\cite{wood},~Lemma 6.2}\textup{)}\label{lemadd}
       Given $A\in \mathbb{R}^{m\times n}$ and $\epsilon\in{(0,1)}$. Suppose $\mathbf{S}$ is a $(1\pm\epsilon)$ $\ell_2$-subspace embedding of $A$. Then we have that
		\begin{equation}\label{eq3}
			(1-\epsilon)\sigma_i(A)\leq \sigma_i(\mathbf{S}A)\leq (1+\epsilon)\sigma_i(A), \text{for all}~ i\in {[n]}.
		\end{equation}
	\end{lem}

	\begin{lem} \textup{(}{\cite{du},Theorem 3.1}\textup{)}\label{lem2}
		Let $x_\star$ denote the precise solution of least squares problems  \eqref{eq1}. For any $x_0 \in R{(A^\ast)}$, the sequence $\{x_k\}_{k=1}^{\infty}$ is  generated by the \textup{MWRK} method. Then
		\begin{equation}\label{eq4}
			\|x_1-x_\star\|_2^2\leq (1-\frac{\sigma_r^2(A)}{\|A\|_F^2})\|x_0-x_\star\|_2^2
		\end{equation}
		and for $k=1, 2, \cdots$,
		\begin{equation}\label{eq5}
			\|x_k-x_\star\|_2^2\leq (1-\frac{\sigma_r^2(A)}{\max\limits_{1\leq j\leq m}\sum\limits_{i=1,i\neq j}^m\|A^{(i)}\|_2^2 })\|x_{k-1}-x_\star\|_2^2
		\end{equation}	
	with the probability at least $1-\delta$.
	\end{lem}

Let  $\{x_k\}_{k=1}^{\infty}$ be generated by the \textup{RaBK} method and $x_k^{\star}=\prod_\chi(x_k)$, where $\chi$ is the solution set of the linear systems $Ax=b$. Suppose that $J_k=J$ for $k=1,2,\cdots$ with the size $|J|=\tau$. In this case, we have $p_i=\frac{\tau}{m}$. Let $\omega_i=\frac{1}{\tau}$ for any $i\in J$, $\|A^{(i)}\|_2=1$ for all $i\in[m]$ and then we have $\|A\|_F^2=m$. For two cases of the stepsize $\alpha_k$ in Remark \ref{alpha}, we give the following lemma.


{
	\begin{lem}\textup{(}{\cite{nec}, Table 1}\textup{)}
		If the sequence $\{\alpha_k\}_{k\geq0}$ is given by case \textup{\textbf{(i)}} or case \textup{\textbf{(ii)}} in Remark \ref{alpha}, we get
		\begin{equation} \label{re1}
			E[\|x_k-x_k^{\star}\|_2^2]\leq
			(1-\frac{\tau}{{\lambda}_{\textup{max}}^{\textup{block}}}\frac{\lambda_{\textup{min}}^{nz}(A^\ast A)}{m})^k\|x_0-x_0^{\star}\|_2^2.
		\end{equation}
		with the probability at least $1-\delta$.
	\end{lem}
}

\section{The proposed methods}\label{section3}
 In this section,  we propose accelerated Kaczmarz methods with randomized sketch techniques for solving consistent linear systems \eqref{eq1} and further consider two types of randomized matrix.

 \subsection{Accelerated Kaczmarz methods with randomized matrix}

Now, based on the frameworks of Algorithm \ref{alg:Framwork1} and Algorithm \ref{alg:Framwork2}, we propose Algorithm \ref{alg:Framwork3} and Algorithm \ref{alg:Framwork4}, respectively.	
 \par
 \begin{algorithm}[htb]
 	\caption{The randomized sketch maximal weighted residual Kaczmarz method(RS-MWRK(S))}
 	\label{alg:Framwork3}
 	\begin{algorithmic}[1]
 		\REQUIRE ~~$ A,b,d~\text{and}~x_0 $\\		
 		\STATE Introduce the randomized matrix $\mathbf{S} \in \mathbb{R}^{d\times m}, d<m$;
 		\STATE Calculate $\tilde{A} = \mathbf{S}A \in {\mathbb{R}}^{d\times n}$ and $\tilde{b}=\mathbf{S}b \in \mathbb{R}^{d}$;
 		\STATE For $k=0, 1, 2, \cdots $ do until satisfy the stopping criteria
 		\STATE Compute $i_k = \arg \max\limits_{1\leq i_k\leq d}\{\frac{\mid \tilde{b}^{(i_k)}-\tilde{A}^{(i_k)}x_k\mid^2}{\parallel \tilde{A}^{(i_k)}\parallel _2 ^2} \}$;
 		\STATE Set $x_{k+1} = x_k+\frac{\tilde{b}^{(i_k)}-\tilde{A}^{(i_k)}x_k}{\parallel \tilde{A}^{(i_k)}\parallel_2^2}(\tilde{A}^{(i_k)})^\ast$;
 		\STATE end for
 			\RETURN $x_{k+1}$;
 	\end{algorithmic}
 \end{algorithm}
 \par
 \begin{algorithm}[htb]
 	\caption{The leverage score sampling randomized average block Kaczmarz method(LS-RaBK(S))}
 	\label{alg:Framwork4}
 	\begin{algorithmic}[1]
 		\REQUIRE ~~	$ A,b,x_0,\{\alpha_k\}_{k\geq0}~\text{and}~\{\omega_k\}_{k\geq0} $\\
 		\STATE Introduce the transform matrix $\mathbf{S} \in \mathbb{R}^{d\times m}, d<m$;
 		\STATE Calculate $\tilde{A} = \mathbf{S}A \in {\mathbb{R}}^{d\times n}$ and $\tilde{b}=\mathbf{S}b \in \mathbb{R}^{d}$;
 		\STATE For $k=0, 1, 2, \cdots $ do until satisfy the stopping criteria
 		\STATE Draw sample $J_k\sim \textup{P}$ and update:
 		\STATE $x_{k+1} = x_k-\alpha_k(\sum\limits_{i\in J_k}\omega_k^i\frac{\tilde{A}^{(i)} x_k-\tilde{b}^{(i)}}{\|\tilde{A}^{(i)}\|_2^2}(\tilde{A}^{(i)})^\ast)$;
 		\STATE end for
 		\RETURN $x_{k+1}$;
 	\end{algorithmic}
 \end{algorithm}

\subsection{Selection of the randomized matrix $\mathbf{S}$}
In the following, two randomized sketch matrices are considered.

\begin{description}
	\item[(I)] Randomized matrix $\mathbf{G}$: A new count sketch matrix is defined as $\mathbf{G} = \mathbf{C} \Phi \in \mathbb{R}^{d\times m}, d<m$, where $\mathbf{C}$ is a $d\times{d}$ diagonal matrix with each diagonal entry independently chosen to be $+1$ or $-1$ with the equal probability and $\Phi \in \{0,1\}^{d\times{m}}$ is a
	$d\times{m}$ binary matrix with $\Phi_{h(i),i} = 1$ and all remaining entries being $0$, where $h:[m]\rightarrow[d]$ is a random map such that for each $i\in [m]$, $h(i)=j$ with probability $\frac{1}{d}$ for each $j \in [d]$.
	\item[(II)] Randomized matrix $\mathbf{Q}$:  A randomized matrix $\mathbf{Q}\in\{0,1\}^{d\times{m}}$ is a binary matrix with $\mathbf{Q}_{i,w(i)}=1$ and all remaining entries being $0$, where $w:[d]\rightarrow[m]$ is a random map such that for each $i\in[d]$, $w(i)=j$ with $j_1\neq j_2\cdots\neq j_d$ for each $j \in [m]$.
\end{description}

\begin{rem}\label{rem}
It is not difficult to check that $\mathbf{G}$ and
$\mathbf{Q}$ satisfy
the conditions of Lemma \ref{lem}.
\end{rem}

For the randomized matrix $\mathbf{G}$, the random $2$-wise independent hash function $h$ determines which row in $\mathbf{G}$ contains the nonzero entry for each column $j\in[m]$.
The randomized sketch matrix $\mathbf{G}$ is different from the one in \cite{zhang}. The product order between the diagonal matrix $\mathbf{C}$ and the binary matrix $\Phi$ is different.  The randomized matrix $\mathbf{G}$ is generated as the product of a $d\times d$ matrix and a $d\times m$ matrix. Due to $d\ll m$, the volume of the $\mathbf{C}$ is less than that of the one in \cite{zhang}, which reduces the calculation time. For large-scale numerical experiments, the proposed methods accelerated by the randomized matrix $\mathbf{G}$ perform better in terms of the running times.

For the randomized matrix $\mathbf{Q}$, we have two remarks: (i) The random $2$-wise independent function $w$ is employed to determine which column of $\mathbf{Q}$ contains the nonzero entry for each row $i\in[d]$. Moreover, the nonzero entries are located in different columns. (ii) Motivated by \cite{wang}, given a matrix $A \in \mathbb{R}^{m\times n}$ with $m\gg n$, and its singular value decomposition $A=UDV^\ast$, the leverage scores of each row are identical. Under such circumstances, instead of selecting rows through leverage scores in \cite{wang}, we can compute the $A_{\mathbf{Q}}=\mathbf{Q}A$ by utilizing the MATLAB function \textit{randperm(m,d)} as follows:

\vspace{1em}
\hrule
\begin{tabular}{l}
     $R=randperm(m,d)$; ~$A_{\mathbf{Q}}=zeros(d,n)$;\\
     for $i=1:d$\\
     ~~~~$A_{\mathbf{Q}}(i,:)=A_{\mathbf{Q}}(i,:)+A(R(i),:);$\\
     end\\
\end{tabular}
\hrule
\vspace{1em}
Furthermore, it is noted that the time complexities of computing $\mathbf{G}A$ and $\mathbf{Q}A$ are $O(nnz(A))$ and $O(nnz(A(d)))$ respectively, where $nnz(A)$ denotes the number of non-zero elements of the matrix $A$ and $nnz(A(d))$ denotes the number of non-zero elements in the selected $d$ rows of matrix $A$.
%
In order to present the theoretical analysis and numerical experiments clearly, we make the following remarks.
 \begin{rem}
	\begin{description}
		\item[(i)] By a large number of numerical experiments, Algorithm \ref{alg:Framwork4} performs worse even does not converge for large-scale experiments. Therefore, we only consider Algorithm \ref{alg:Framwork4} in the case that $\mathbf{S}=\mathbf{Q}$.
		\item[(ii)] For convenience, if $\mathbf{S}=\mathbf{G}$(or $\mathbf{S}=\mathbf{Q}$), \textup{RS-MWRK(S)} is renamed as \textup{RS-MWRK(G)}(or \textup{RS-MWRK(Q)}). For Algorithm \ref{alg:Framwork4}, \textup{LS-RaBK(S)} is renamed  as \textup{LS-RaBK(Q)} if $\mathbf{S}=\mathbf{Q}$.	
	\end{description}
\end{rem}

\section{Convergence analysis}\label{section4}
	In this section, we give the convergence analysis of the proposed methods.
{
\begin{lem} \label{lem5}
Let $\mathbf{S}=\mathbf{G}$(or $\mathbf{S}=\mathbf{Q}$) with $d=O(n^2/(\delta\epsilon^2))$, where $\epsilon,\delta\in{(0,1)}$. For any $A \in \mathbb{R}^{m\times n}$, we have $R(A^\ast\mathbf{S}^\ast)=R(A^\ast)$ with probability $1-\delta$.
\end{lem}}
\begin{proof} It is not difficult to check that $R(A^\ast\mathbf{S}^\ast)\subseteq R(A^\ast)$. We just prove $R(A^\ast\mathbf{S}^\ast)^\perp\subseteq R(A^\ast)^\perp$. For any $x\in R(A^\ast\mathbf{S}^\ast)^\perp$, we have $\mathbf{S}Ax = 0$, which yields $\|\mathbf{S}Ax\|_2 = 0$. By Lemma $\ref{lem}$, Lemma $\ref{lem1}$ and Remark $\ref{rem}$, one can check that $\|Ax\|_2 = 0$ with probability $1-\delta$. Thus, we derive $x\in R(A^\ast)^\perp$ which implies that $R(A^\ast\mathbf{S}^\ast)^\perp\subseteq R(A^\ast)^\perp$ with probability $1-\delta$.
\end{proof}

\begin{thm}\label{thm}
Let $\mathbf{S}=\mathbf{G}$(or $\mathbf{S}=\mathbf{Q}$) with $d=O(n^2/(\delta\epsilon^2))$, where $\epsilon,\delta\in{(0,1)}$, $x_\star$ denote the precise solution of least squares problems  \eqref{eq1}, $\{x_k\}_{k=1}^{\infty}$ with $x_0 \in R{(A^\ast)}$ denote the iteration sequence generated by the \textup{RS-MWRK(S)} method and $\tilde{A} = \mathbf{S}A$. Then we have 
	\begin{equation}\label{eq6}
		\|x_1-x_\star\|_2^2\leq (1-\frac{(1-\epsilon)^2}{(1+\epsilon)^2}\frac{\sigma_r^2(A)}{\|A\|_F^2})\|x_0-x_\star\|_2^2,
	\end{equation}
	and for $k=1, 2, \cdots$,
	\begin{equation}\label{eq7}
		\|x_k-x_\star\|_2^2\leq (1-(1-\epsilon)^2\frac{\sigma_r^2(A)}{\max\limits_{1\leq j\leq d}\sum\limits_{i=1,i\neq j}^d\|\tilde{A}^{(i)}\|_2^2 })\|x_{k-1}-x_\star\|_2^2
	\end{equation}
	 with probability at least $1-\delta$.
\end{thm}
	\begin{proof}
By	Lemmas \ref{lem}-\ref{lem1}, Lemma \ref{lem2} and Lemma \ref{lem5}, the desired assertion can be given via employing the similar proof technique of Theorem 3 in the reference \cite{zhang}.		
	\end{proof}

Next we give an upper bound of the distance between the numerical solution derived from the \textup{MWRK} method and the one from the \textup{RS-MWRK} method.
	\begin{thm}\label{thm2}
		Let $\mathbf{S}=\mathbf{G}$  (or $\mathbf{S}=\mathbf{Q}$) with $d=O(n^2/(\delta\epsilon^2))$, where $\epsilon,\delta\in{(0,1)}$ and $x_\star$ denote the precise solution of least squares problems  \eqref{eq1}. The sequence $\{x_k\}_{k=1}^{\infty}$ is generated by the \textup{MWRK} method and $\{x_k^{'}\}_{k=1}^{\infty}$ denotes the iteration sequence generated by the \textup{RS-MWRK(S)} method. Then we get
		\begin{equation}\label{eq12}
			\|x_k^{'}-x_k\|_2^2\leq (4-(4-6\epsilon+3\epsilon^2)\frac{\sigma_r^2(A)}{t})p,
		\end{equation}
		 with probability at least $1-\delta$, where $t=\max\{\max\limits_{1\leq j\leq m}\sum\limits_{i=1,i\neq j}^m\|A^{(i)}\|_2^2, \max\limits_{1\leq j\leq d}\sum\limits_{i=1,i\neq j}^d\|\tilde{A}^{(i)}\|_2^2\}$ and $p=\max\{\|x_{k-1}^{'}-x_\star\|_2^2,\|x_{k-1}-x_\star\|_2^2\}$.
	\end{thm}	
	\begin{proof} Note that
{
     \begin{equation}\label{eq24}
     \|x_k^{'}-x_k\|_2^2=\|x_k^{'}-x_\star+x_\star-x_k\|_2^2=\|x_k^{'}-x_\star\|_2^2 + 2\langle x_k^{'}-x_\star,x_\star - x_k\rangle+\|x_\star - x_k\|_2^2.
     \end{equation}}
     \par
     By the Cauchy-Schwarz inequality, Lemma \ref{lem2} and Theorem \ref{thm}, with probability at least $1-\delta$, we have
     $$\langle x_k^{'}-x_\star,x_\star - x_k\rangle\leq \|x_k^{'}-x_\star\|_2 \|x_\star - x_k\|_2\leq \sqrt{(1-(1-\epsilon)^2\frac{\sigma_r^2(A)}{t})(1-\frac{\sigma_r^2(A)}{t})}~p.$$
     Due to $1-(1-\epsilon)^2\frac{\sigma_r^2(A)}{t}\geq 1-\frac{\sigma_r^2(A)}{t}$, it is derived that
     \begin{equation}\label{eq25}
     \langle x_k^{'}-x_\star,x_\star - x_k\rangle\leq (1-(1-\epsilon)^2\frac{\sigma_r^2(A)}{t})p.
     \end{equation}
     From the equalities \eqref{eq5} and \eqref{eq7}, with probability at least $1-\delta$, one can check that
     \begin{equation}\label{eq26}
     \|x_k^{'}-x_\star\|_2^2+\|x_\star - x_k\|_2^2\leq (2-(2-2\epsilon+\epsilon^2)\frac{\sigma_r^2(A)}{t})p,
     \end{equation}
    together with \eqref{eq25} and \eqref{eq24} yields
    $$\|x_k^{'}-x_k\|_2^2\leq (4-(4-6\epsilon+3\epsilon^2)\frac{\sigma_r^2(A)}{t})p$$
     with probability at least $1-\delta$.
\end{proof}

 Let $\mathbf{S}=\mathbf{Q}$ with $d=O(n^2/(\delta\epsilon^2))$ where $\epsilon,\delta\in{(0,1)}$. Let $\{x_k\}_{k=1}^{\infty}$ denote the iteration sequence generated by the \textup{LS-RaBK(Q)} method and $x_k^{\star}=\prod_\chi(x_k)$ with $x_0 \in R{(A^\ast)}$, {where $\chi$ is the solution set of the linear systems $Ax=b$.} Suppose that $J_k=J$ for $k=1,2,\cdots $ with the size $|J|={\tau}$.  In this case, we have $p_i=\frac{\tau}{d}$. Set $\omega_i=\frac{1}{\tau}$. Next we present the convergence analysis of the LS-RaBK(Q) method based on two selection methods of parameter $\alpha_k$.

	\begin{thm}\label{thm1}
		
		If the sequence $\{\alpha_k\}_{k\geq0}$ is given by case \textup{\textbf{(i)}} or case \textup{\textbf{(ii)}} in Remark \ref{alpha}, then $\{x_k\}_{k=1}^{\infty}$ generated by the \textup{LS-RaBK(Q)} method satisfies	
			\begin{equation}\label{re2}
		E[\|x_k-x_k^{\star}\|_2^2]\leq
		(1-\frac{(1-\epsilon)^2}{(1+\epsilon)^2}\frac{\tau}{\lambda_{\textup{max}}^{\textup{block}}}\frac{\lambda_{\textup{min}}^{nz}(A^\ast A)}{d})^k\|x_0-x_0^{\star}\|_2^2
	\end{equation}
		with probability at least $1-\delta$.
	\end{thm}
	\begin{proof}
		By Lemma \ref{lem5}, replacing the matrix $A$ by $\tilde{A}$ in the inequality \eqref{re1}, with probability at least $1-\delta$, it is derived that
		\begin{equation}\label{eq16}
			E[\|x_k-x_k^{\star}\|_2^2]\leq
			(1-\frac{\tau}{\tilde{\lambda}_{\textup{max}}^{\textup{block}}}\frac{\lambda_{\textup{min}}^{nz}(A^\ast \mathbf{S}^\ast \mathbf{S} A)}{d})^k\|x_0-x_0^{\star}\|_2^2,
		\end{equation}
	where $\tilde{\lambda}_{\textup{max}}^{\textup{block}}=\max\limits_{J\sim \textup{P}}\lambda_{\textup{max}}(\tilde{A}_J^\ast \tilde{A}_J)$. From Lemma \ref{lemadd}, with probability at least $1-\delta$, one can check that
		
		\begin{equation}\label{eq17}
			\tilde{\lambda}_{\textup{max}}^{\textup{block}}=\max\limits_{J\sim \textup{P}}\lambda_{\textup{max}}(\tilde{A_J^\ast} \tilde{A_J})=\max\limits_{J\sim \textup{P}}\sigma_{\textup{max}}^2(\tilde{A_J})
			\leq\max\limits_{J\sim \textup{P}}(1+\epsilon)^2\sigma_{\textup{max}}^2(A_J)=(1+\epsilon)^2\lambda_{\textup{max}}^{\textup{block}},
		\end{equation}
	and
		\begin{equation}\label{eq18}
			\lambda_{\textup{min}}^{nz}(A^\ast \mathbf{S}^\ast \mathbf{S} A)=(\sigma_{\textup{min}}^{{nz}}(\mathbf{S}A))^2\geq(1-\epsilon)^2(\sigma_{\textup{min}}^{{nz}}(A))^2=(1-\epsilon)^2\lambda_{\textup{min}}^{nz}(A^\ast A)
		\end{equation}
		which together with the inequality \eqref{eq16} derives the inequality \eqref{re2} with probability at least $1-\delta$.
	\end{proof}

Furthermore, we can also provide upper bounds for the distance between the numerical solution derived from the \textup{RaBk} method and that from the \textup{LS-RaBK(Q)} method, considering different selections of the parameter $\alpha_k$.
{
		\begin{thm}
		Let $\mathbf{S}=\mathbf{Q}$ with $d=O(n^2/(\delta\epsilon^2))$ where $\epsilon,\delta\in{(0,1)}$. The sequence $\{x_k\}_{k=1}^{\infty}$ is generated by the \textup{RaBK} method and $\{x_k^{'}\}_{k=1}^{\infty}$ denotes the iteration sequence generated by the \textup{LS-RaBK(Q)} method. If the sequence $\{\alpha_k\}_{k\geq0}$ is given by case \textup{\textbf{(i)}} or case \textup{\textbf{(ii)}} in Remark \ref{alpha}, then we yield
			\begin{equation}\label{re3}
				E[\|x_k^{'}-x_k\|_2^2]\leq ((1-\frac{\tau}{\lambda_{\textup{max}}^{\textup{block}}}\frac{\lambda_{\textup{min}}
					^{nz}(A^\ast A)}{m})^k+3(1-\frac{(1-\epsilon)^2}{(1+\epsilon)^2}\frac{\tau}{\lambda_{\textup{max}}^{\textup{block}}}\frac{\lambda_{\textup{min}}^{nz}(A^\ast A)}{m})^k)\|x_0-x_0^{\star}\|_2^2
			\end{equation}
			with probability at least $1-\delta$.
	\end{thm}
}

\begin{proof} Since $d<m$, we have $$1-\frac{(1-\epsilon)^2}{(1+\epsilon)^2}\frac{\tau}{\lambda_{\textup{max}}^{\textup{block}}}\frac{\lambda_{\textup{min}}^{nz}(A^\ast A)}{d}<1-\frac{(1-\epsilon)^2}{(1+\epsilon)^2}\frac{\tau}{\lambda_{\textup{max}}^{\textup{block}}}\frac{\lambda_{\textup{min}}^{nz}(A^\ast A)}{m},$$
combing with
$$1-\frac{\tau}{\lambda_{\textup{max}}^{\textup{block}}}\frac{\lambda_{\textup{min}}
	^{nz}(A^\ast A)}{m}\le1-\frac{(1-\epsilon)^2}{(1+\epsilon)^2}\frac{\tau}{\lambda_{\textup{max}}^{\textup{block}}}\frac{\lambda_{\textup{min}}^{nz}(A^\ast A)}{m}$$	
gives the inequality \eqref{re3} by a similar proof with Theorem \ref{thm2}.
\end{proof}
\section{Numerical experiments}\label{section5}
All the experiments in this section are performed in MATLAB R2022a on a personal computer with 2.50GHz central processing unit (Intel(R) Core(TM) i5-12500H CPU), 16.00 GB memory, and Windows operating system (Windows 11). Two aspects are given to check the efficiency of the proposed algorithms: the number
of iteration steps (IT), the running time in seconds (CPU). It is noted that `IT' and `CPU' are respectively the average number of the iteration steps and the running time when the tested method runs fifty times. To facilitate the comparison of the running time for the tested methods, we define
$$ \mbox{CPU speedup 1}=\frac{\mbox{CPU of CS-MWRK}}{\mbox{CPU of RS-MWRK(Q)}}, $$
$$ \mbox{CPU speedup 2}=\frac{\mbox{CPU of CS-MWRK}}{\mbox{CPU of RS-MWRK(G)}}, $$
$$ \mbox{CPU speedup 3}=\frac{\mbox{CPU of RaBK-c}}{\mbox{CPU of LS-RaBK(Q)-c}}, $$
and
$$ \mbox{CPU speedup 4}=\frac{\mbox{CPU of RaBK-a}}{\mbox{CPU of LS-RaBK(Q)-a}}. $$
\par
In the following, the right-hand side vector $b=Ax_\star$, where $A$ and $x_\star$ are generated by the MATLAB function \textit{randn}. For the uniform sampling variants, we choose $\tau=\frac{d}{50}$ and all the experiments start from the initial vector $x_0=0$. Every test algorithm is terminated once the relative solution error (RES) at the current iteration satisfies $\mbox{RES}<10^{-6}$ or when the number of iteration steps exceeds $100000$, where

$$ \mbox{RES} = \frac{\|x_k-x_\star\|_2^2}{\|x_\star\|_2^2}. $$
\par
The test algorithms are listed in Table \ref{add}.
\begin{table}[h]
	\caption{Abbreviations of the tested methods}
	\centering
	\scalebox{0.75}{
	\begin{tabular}{lll}
		\hline
		Methods &$\alpha_k$ &Abbreviations\\
		\hline
		Randomized sketch maximal weighted residual Kaczmarz method (S=G)&-&RS-MWRK(G)\\
		Randomized sketch maximal weighted residual Kaczmarz method (S=Q)&-&RS-MWRK(Q)\\
		Leverage score sampling randomized average block Kaczmarz method with the constant stepsize &$1.95$&LS-RaBK(Q)-c\\
		Leverage score sampling randomized average block Kaczmarz method with the adaptive stepsize &$1.95L_k$&LS-RaBK(Q)-a\\
		Count sketch maximal weighted residual Kaczmarz method in \cite{zhang}&-&CS-MWRK\\
		Randomized average block Kaczmarz method with the constant stepsize in \cite{nec}&$1.95$&RaBK-c\\
		Randomized average block Kaczmarz method with the adaptive stepsize in \cite{nec}&$1.95L_k$&RaBK-a\\
		\hline
	\end{tabular}}
	\label{add}
\end{table}

\begin{sidewaystable}[htbp]
	\centering
	\caption{Numerical results of CS-MWRK, RS-MWRK(Q), RS-MWRK(G), LS-RaBK(Q)-c, RaBK-c, LS-RaBK(Q)-a and RaBK-a for different random matrices $A$.}
	\scalebox{1}{
		\begin{tabular}{cccccccccc}
			\hline
			$m\times n$ &$d$ & &CS-MWRK &RS-MWRK(Q) &RS-MWRK(G) &LS-RaBK(Q)-c &RaBK-c &LS-RaBK(Q)-a &RaBK-a  \\
			\hline
			$5000\times 50$ &$10n$ &IT &85.3600 &85.9400 &85.9000 &281.3200 &203.8600 &1536.3200 &1119.6400  \\
			&&CPU &0.0025 &\textbf{0.0006} &0.0019 &0.0016 &0.0025 &0.0159 &0.0072\\
			\cmidrule(lr){2-10}
			&$20n$ &IT &67.8200 &67.5200 &68.3200 &229.3200 &191.2200 &725.2800 &638.2000 \\
			&&CPU &0.0025 &0.0050 &\textbf{0.0016} &0.0022 &0.0034 &0.0122 &0.0059\\
			\cmidrule(lr){2-10}
			&$30n$ &IT &61.4000 &60.9400 &61.0800 &207.5600 &185.6400 &511.5800 &473.7200  \\
			&&CPU &0.0047 &0.0084 &0.0028 &0.0044 &\textbf{0.0013} &0.0225 &0.0147\\
			\cmidrule(lr){2-10}
			&$40n$ &IT &58.3800 &57.3800 &58.5000 &202.7600 &187.2600 &432.6200 &405.1400 \\
			&&CPU &0.0016 &0.0028 &\textbf{0.0006} &0.0022 &0.0019 &0.0156 &0.0153\\
			\hline
			$5000\times 100$ &$5n$ &IT &264.2600 &267.4200 &268.6200 &803.5000 &423.4400 &8098.0200 &4281.9600 \\
			&&CPU &0.0044 &0.0191 &\textbf{0.0031} &0.0050 &0.0038 &0.1134 &0.0675\\
			\cmidrule(lr){2-10}
			&$10n$ &IT &173.7200 &172.4000 &172.4000 &577.0400 &421.6400 &3116.3800 &2399.5200 \\
			&&CPU &\textbf{0.0019} &0.0100 &0.0028 &0.0075 &0.0044 &0.0853 &0.0697\\
			\cmidrule(lr){2-10}
			&$15n$ &IT & 148.8800  & 145.6000  & 148.0400  & 487.0800  & 402.4400  & 1962.1800  & 1671.8400 \\
			&&CPU & 0.0022  & \textbf{0.0019}  & 0.0025  & 0.0094  & 0.0053  & 0.0825  & 0.0650\\
			\cmidrule(lr){2-10}
			&$20n$ &IT &136.5400  & 134.0800  & 137.1600  & 461.8000  & 404.0000  & 1483.0800  & 1347.1200 \\
			&&CPU & 0.0028  & \textbf{0.0025}  & 0.0028  & 0.0194  & 0.0100  & 0.0834  & 0.0775\\
			\hline
			$5000\times 150$ &$5n$ &IT & 393.8200  & 400.9800  & 396.1200  & 1205.6000  & 630.9000  & 12285.9800  & 6542.1600 \\
			&&CPU & \textbf{0.0038}  & 0.0075  & \textbf{0.0038}  & 0.0191  & 0.0059  & 0.3016  & 0.1819\\
			\cmidrule(lr){2-10}
			&$10n$ &IT &261.4400  & 255.2600  & 260.8200  & 830.5000  & 623.2800  & 4552.2000  & 3558.0800 \\
			&&CPU & \textbf{0.0041}  & 0.0069  & 0.0063  & 0.0219  & 0.0181  & 0.2419  & 0.1875\\
			\cmidrule(lr){2-10}
			&$15n$ &IT &228.3600  & 218.2400  & 226.6000  & 733.9400  & 631.0000  & 2936.8800  & 2595.9600 \\
			&&CPU & 0.0045  & 0.0059  & \textbf{0,0028} & 0.0378  & 0.0322  & 0.2250  & 0.2069\\
			\cmidrule(lr){2-10}
			&$20n$ &IT &214.5600  & 199.1800  & 213.0000  & 698.2000  & 635.8200  & 2239.1600  & 2092.4400 \\
			&&CPU & \textbf{0.0103}  & \textbf{0.0103}  & 0.0106  & 0.0500  & 0.0328  & 0.2366  & 0.2163 \\
			\hline
	\end{tabular}}
	\label{tab1}
\end{sidewaystable}

\begin{table}[htbp]
     \centering
     \caption{Numerical results of CPU speedup 1, CPU speedup 2, CPU speedup 3 and CPU speedup 4 for different scale matrices $A$.}
     \begin{tabular}{cccccc}
     \hline
     $m\times n$ &$d$ &CPU speedup 1 &CPU speedup 2 &CPU speedup 3 &CPU speedup 4  \\
     \hline
     $5000\times 50$ &$10n$ &\textbf{4.1667}  & 1.3158  & \textbf{1.5625}  & 0.4528  \\
     &$20n$ &0.5000  & 1.5625  & 1.5455  & 0.4836  \\
     &$30n$ &0.5595  & 1.6786  & 0.2955  & 0.6533  \\
     &$40n$ &0.5714  & \textbf{2.6667}  & 0.8636  & \textbf{0.9808}  \\
     \hline
     $5000\times 100$ &$5n$ & 0.2304  & 1.4194  & 0.7600  & 0.5952  \\
     &$10n$ &0.1900  & 0.6786  & 0.5867  & 0.8171  \\
     &$15n$ &1.1579  & 0.8800  & 0.5638  & 0.7879  \\
     &$20n$ &1.1200  & 1.0000  & 0.5155  & 0.9293  \\
     \hline
     $5000\times 150$ &$5n$ &0.5067  & 1.0000  & 0.3089  & 0.6031  \\
     &$10n$ &0.5942  & 0.6508  & 0.8265  & 0.7751  \\
     &$15n$ &0.7627  & 1.6071  & 0.8519  & 0.9196  \\
     &$20n$ &1.0000  & 0.9717  & 0.6560  & 0.9142  \\
     \hline
     \end{tabular}
 \label{addt1}
\end{table}

\begin{sidewaystable}[htbp]
	\centering
	\caption{Numerical results of CS-MWRK, RS-MWRK(Q), RS-MWRK(G), LS-RaBK(Q)-c, RaBK-c, LS-RaBK(Q)-a and RaBK-a for different  scale
		 matrices $A$.}
	\scalebox{1}{
		\begin{tabular}{cccccccccc}
			\hline
			$m\times n$ &$d$ & &CS-MWRK &RS-MWRK(Q) &RS-MWRK(G) &LS-RaBK(Q)-c &RaBK-c &LS-RaBK(Q)-a &RaBK-a  \\
			\hline
			$50000\times 50$ &$20n$ &IT & 68.1600  & 67.4200  & 67.0000  & 229.1800  & 184.7800  & 734.9400  & 616.1200  \\
			&&CPU & 0.0128  & 0.0113  & 0.0088  & \textbf{0.0038}  & 0.0141  & 0.0238  & 0.0331\\
			\cmidrule(lr){2-10}
			&$30n$ &IT & 60.5400  & 61.0000  & 61.2600  & 212.8400  & 182.8200  & 520.2200  & 469.7800 \\
			&&CPU & 0.0163  & \textbf{0.0028}  & 0.0153  & 0.0034  & 0.0119  & 0.0175  & 0.0281\\
			\cmidrule(lr){2-10}
			&$40n$ &IT & 57.9600  & 57.0600  & 57.4600  & 203.6600  & 179.3400  & 425.2600  & 395.5800  \\
			&&CPU & 0.0147  & \textbf{0.0038}  & 0.0088  & 0.0100  & 0.0222  & 0.0175  & 0.0341\\
			\cmidrule(lr){2-10}
			&$50n$ &IT & 55.3800  & 54.7000  & 54.7400  & 198.6800  & 179.5600  & 371.8800  & 351.6200  \\
			&&CPU & 0.0128  & \textbf{0.0038}  & 0.0100  & 0.0097  & 0.0294  & 0.0169  & 0.0369\\
			\hline
			$50000\times 100$ &$10n$ &IT & 172.1200  & 170.2600  & 172.0800  & 563.2800  & 373.6200  & 3057.7400  & 2138.6000  \\
			&&CPU & 0.0303  & \textbf{0.0066}  & 0.0200  & 0.0138  & 0.0453  & 0.0800  & 0.1234\\
			\cmidrule(lr){2-10}
			&$20n$ &IT & 133.5000  & 133.3800  & 135.2600  & 462.1400  & 368.2400  & 1485.0200  & 1254.0600  \\
			&&CPU & 0.0297  & \textbf{0.0116}  & 0.0231  & 0.0225  & 0.0434  & 0.0850  & 0.1113\\
			\cmidrule(lr){2-10}
			&$50n$ &IT & 107.6400  & 107.7400  & 108.5600  & 399.4800  & 360.9600  & 748.7200  & 711.8800  \\
			&&CPU & 0.0488  & \textbf{0.0175}  & 0.0303  & 0.0369  & 0.0653  & 0.1106  & 0.1306\\
			\cmidrule(lr){2-10}
			&$70n$ &IT & 101.1400  & 101.0600  & 101.0400  & 388.2000  & 359.7400  & 625.3600  & 607.0200 \\
			&&CPU & 0.0478  & \textbf{0.0225}  & 0.0388  & 0.0622  & 0.0903  & 0.1494  & 0.1569\\
			\cmidrule(lr){2-10}
			&$100n$ &IT & 95.4600  & 95.4600  & 95.3400  & 376.0400  & 358.5200  & 539.7000  & 530.8600 \\
			&&CPU & 0.0703  & \textbf{0.0447}  & 0.0603  & 0.0919  & 0.1047  & 0.1925  & 0.2116\\
			\hline
			$50000\times 150$ &$10n$ &IT & 255.3800  & 255.3800  & 256.6800  & 838.5600  & 553.6400  & 4616.1800  & 3208.2800 \\
			&&CPU & 0.0431  & \textbf{0.0116}  & 0.0394  & 0.0281  & 0.0819  & 0.2450  & 0.2928\\
			\cmidrule(lr){2-10}
			&$50n$ &IT & 161.2800  & 160.7000  & 160.6400  & 598.5200  & 545.5000  & 1124.7000  & 1073.3800 \\
			&&CPU & 0.1163  & \textbf{0.0566}  & 0.0750  & 0.1141  & 0.1622  & 0.2988  & 0.3678\\
			\cmidrule(lr){2-10}
			&$100n$ &IT & 142.6000  & 141.2000  & 142.4600  & 565.0000  & 540.2600  & 809.7600  & 798.7000 \\
			&&CPU & 0.0919  & \textbf{0.0669}  & 0.0916  & 0.2453  & 0.2681  & 0.5022  & 0.5338\\
			\cmidrule(lr){2-10}
			&$150n$ &IT & 135.0000  & 133.2600  & 135.8600  & 559.6800  & 548.3000  & 715.4000  & 711.7200 \\
			&&CPU & 0.2000  & \textbf{0.1363}  & 0.1731  & 0.3916  & 0.3994  & 0.6919  & 0.6931 \\
			\hline
	\end{tabular}}
	\label{tab2}
\end{sidewaystable}

\begin{table}[htbp]
     \centering
     \caption{Numerical results of CPU speedup 1, CPU speedup 2, CPU speedup 3 and CPU speedup 4 for different scale matrices $A$.}
     \begin{tabular}{cccccc}
     \hline
     $m\times n$ &$d$ &CPU speedup 1 &CPU speedup 2 &CPU speedup 3 &CPU speedup 4  \\
     \hline
     $50000\times 50$ &$20n$ &1.1389  & 1.4643  & \textbf{3.7500}  & 1.3947  \\
     &$30n$ &\textbf{5.7778}  & 1.0612  & 3.4545  & 1.6071  \\
     &$40n$ &3.9167  & \textbf{1.6786}  & 2.2188  & 1.9464  \\
     &$50n$ &3.4167  & 1.2813  & 3.0323  & \textbf{2.1852}  \\
     \hline
     $50000\times 100$ &$10n$ & 4.6190  & 1.5156  & 3.2955  & 1.5430  \\
     &$20n$ &2.5676  & 1.2838  & 1.9306  & 1.3088  \\
     &$50n$ &2.7857  & 1.6082  & 1.7712  & 1.1808  \\
     &$70n$ &2.1250  & 1.2339  & 1.4523  & 1.0502  \\
     &$100n$ &1.5734  & 1.1658  & 1.1395  & 1.0990  \\
     \hline
     $50000\times 150$ &$10n$ &3.7297  & 1.0952  & 2.9111  & 1.1952  \\
     &$50n$ &2.0552  & 1.5500  & 1.4219  & 1.2312  \\
     &$100n$ &1.3738  & 1.0034  & 1.0930  & 1.0629  \\
     &$150n$ &1.4679  & 1.1552  & 1.0200  & 1.0018  \\
     \hline
     \end{tabular}
 \label{addt2}
\end{table}

\begin{sidewaystable}[htbp]
	\centering
	\caption{Numerical results of CS-MWRK, RS-MWRK(Q), RS-MWRK(G), LS-RaBK(Q)-c, RaBK-c, LS-RaBK(Q)-a and RaBK-a  different  scale
		matrices $A$.}
	\scalebox{1}{
		\begin{tabular}{cccccccccc}
			\hline
			$m\times n$ &$d$ & &CS-MWRK &RS-MWRK(Q) &RS-MWRK(G) &LS-RaBK(Q)-c &RaBK-c &LS-RaBK(Q)-a &RaBK-a  \\
			\hline
			$500000\times 50$ &$10n$ &IT & 86.2000  & 86.1000  & 84.9800  & 282.1400  & 192.6000  & 1538.9800  & 1046.7600  \\
			&&CPU & 0.0594  & 0.0078  & 0.0372  & \textbf{0.0022}  & 0.1406  & 0.0178  & 0.1866\\
			\cmidrule(lr){2-10}
			&$20n$ &IT & 67.7000  & 67.1800  & 67.6400  & 228.9400  & 183.6000  & 726.1800  & 612.7000 \\
			&&CPU & 0.0619  & 0.0093  & 0.0466  & \textbf{0.0041}  & 0.1328  & 0.0141  & 0.1622\\
			\cmidrule(lr){2-10}
			&$40n$ &IT & 57.6200  & 56.8000  & 57.1600  & 203.1400  & 179.0800  & 421.7000  & 393.5800  \\
			&&CPU & 0.0713  & \textbf{0.0034}  & 0.0422  & \textbf{0.0034}  & 0.1834  & 0.0169  & 0.1878\\
			\cmidrule(lr){2-10}
			&$50n$ &IT & 54.9000  & 55.2400  & 54.6600  & 198.2800  & 177.7000  & 370.1200  & 349.8600  \\
			&&CPU & 0.0634  & 0.0038  & 0.0547  & \textbf{0.0028}  & 0.1666  & 0.0166  & 0.2053\\
			\hline
			$500000\times 100$ &$20n$ &IT & 134.1400  & 133.3000  & 134.4200  & 455.2800  & 358.9400  & 1473.3400  & 1231.8400  \\
			&&CPU & 0.3516  & \textbf{0.0059}  & 0.1600  & 0.0172  & 0.7141  & 0.0834  & 0.9181\\
			\cmidrule(lr){2-10}
			&$50n$ &IT  & 107.9000  & 108.0600  & 107.2000  & 399.2400  & 355.3800  & 745.0800  & 706.1200  \\
			&&CPU & 0.4159  & \textbf{0.0075}  & 0.1947  & 0.0328  & 0.8281  & 0.1041  & 0.9503\\
			\cmidrule(lr){2-10}
			&$80n$ &IT & 99.1000  & 99.0200  & 98.7000  & 380.7400  & 353.7800  & 588.6400  & 572.2200  \\
			&&CPU & 0.4747  & \textbf{0.0184}  & 0.2591  & 0.0666  & 0.8853  & 0.1456  & 0.9713\\
			\cmidrule(lr){2-10}
			&$100n$ &IT & 94.4400  & 95.2400  & 94.9000  & 374.6600  & 352.5000  & 538.4800  & 527.4000 \\
			&&CPU & 0.4750  & \textbf{0.0256}  & 0.2653  & 0.0819  & 0.9213  & 0.1641  & 1.0106\\
			\hline
			$500000\times 150$ &$20n$ &IT & 200.2600  & 200.3400  & 201.2800  & 689.4800  & 538.0800  & 2233.6200  & 1857.7800 \\
			&&CPU & 0.5516  & \textbf{0.0434}  & 0.3328  & 0.0513  & 1.5484  & 0.3041  & 1.7469 \\
			\cmidrule(lr){2-10}
			&$50n$ &IT & 160.4800  & 160.6200  & 161.0600  & 598.5000  & 531.6200  & 1120.4200  & 1058.8200 \\
			&&CPU & 0.7100  & \textbf{0.0291}  & 0.3731  & 0.1163  & 1.6922  & 0.3031  & 1.9922\\
			\cmidrule(lr){2-10}
			&$100n$ &IT & 141.0000  & 141.9800  & 141.7000  & 565.0000  & 532.2200  & 809.9600  & 795.1400 \\
			&&CPU & 0.6906  & \textbf{0.0731}  & 0.4791  & 0.2619  & 1.8563  & 0.4809  & 2.1213\\
			\cmidrule(lr){2-10}
			&$150n$ &IT & 132.6800  & 133.3400  & 132.4200  & 552.2400  & 530.7600  & 712.7000  & 706.0200 \\
			&&CPU & 0.8669  & \textbf{0.2116}  & 0.6684  & 0.4041  & 1.9778  & 0.7122  & 2.3891 \\
			\hline
	\end{tabular}}
	\label{tab3}
\end{sidewaystable}

\begin{table}[htbp]
     \centering
     \caption{Numerical results of CPU speedup 1, CPU speedup 2, CPU speedup 3 and CPU speedup 4 for different scale matrices $A$.}
     \begin{tabular}{cccccc}
     \hline
     $m\times n$ &$d$ &CPU speedup 1 &CPU speedup 2 &CPU speedup 3 &CPU speedup 4  \\
     \hline
     $500000\times 50$ &$10n$ &7.6000  & 1.5966  & \textbf{64.2857}  & 10.4737  \\
     &$20n$ &6.6532  & 1.3289  & 32.6923  & 11.5333  \\
     &$40n$ &20.9559  & 1.6889  & 53.3636  & 11.1296  \\
     &$50n$ &16.9167  & 1.1600  & 59.2222  & \textbf{12.3962}  \\
     \hline
     $500000\times 100$ &$20n$ & \textbf{59.2105}  & \textbf{2.1973}  & 41.5455  & 11.0037  \\
     &$50n$ &55.4583  & 2.1364  & 25.2381  & 9.1321  \\
     &$80n$ &25.7458  & 1.8323  & 13.3005  & 6.6695  \\
     &$100n$ &18.5366  & 1.7903  & 11.2519  & 6.1600  \\
     \hline
     $500000\times 150$ &$20n$ &12.6978  & 1.6573  & 30.2134  & 5.7451  \\
     &$50n$ &24.4301  & 1.9028  & 14.5565  & 6.5722  \\
     &$100n$ &9.4444  & 1.4416  & 7.0883  & 4.4107  \\
     &$150n$ &4.0975  & 1.2969  & 4.8948  & 3.3545  \\
     \hline
     \end{tabular}
 \label{addtable}
\end{table}

Let $m=5000$, $50000$ and $500000$ respectively. For given $m$, set $n=50,100,150$ respectively. We record the numerical results of the tested algorithms in
Table \ref{tab1}, Table \ref{tab2} and Table \ref{tab3}. Furthermore, we list the  value of CPU speedup 1-CPU speedup 4 in Table \ref{addt1}, Table \ref{addt2} and Table \ref{addtable} that relate to the different scale matrices. It can be seen that the algorithms with the randomized matrix are improved in terms of the running time compared with the original ones.  As the scale of the matrix increases, the runtime of the RS-MWRK(Q) method and the RS-MWRK(G) method are less than that of the CS-MWRK method, and the CPU speedup reaches 59.2105 and 2.6667 in the experiments, respectively. In addition, it is seen that the LS-RaBK(Q)-c method and LS-RaBK(Q)-a method require more iterations than the RaBK-c method and the RaBK-a method. However, the runtime of the LS-RaBK(Q)-c method and LS-RaBK(Q)-a is less than that of the RaBK-c method and the RaBK-a method, and the CPU speedup is 64.2857 and 12.3962, respectively. This is because the LS-RaBK(Q)-c method and LS-RaBK(Q)-a reduce the dimension by means of the randomized sketch matrix. As a result, a smaller matrix is preprocessed and the calculation time is reduced.
\par
In Tables \ref{tab1}-\ref{addt1}, one can check that the proposed algorithm outperforms the other algorithms in most cases in terms of the running time.
From Table \ref{tab2} and Table \ref{tab3}, the runtime of the RS-MWRK(Q) method and RS-MWRK(G) method is lower than that of the CS-MWRK method for a larger $m$. 
 As the size of the matrix increases, the RS-MWRK(G) method is faster than the CS-MWRK method. For the RS-MWRK(Q) method, we only need to randomly select the number of matrix rows without matrix multiplication, which greatly reduces the running time.  We can also see from the numerical experiments that the RS-MWRK(Q) method seems to be more efficient for solving super large scale highly overdetermined linear systems.

\begin{figure}[htbp]
	\centering
	\includegraphics[width=0.4\linewidth,height=0.4\linewidth]{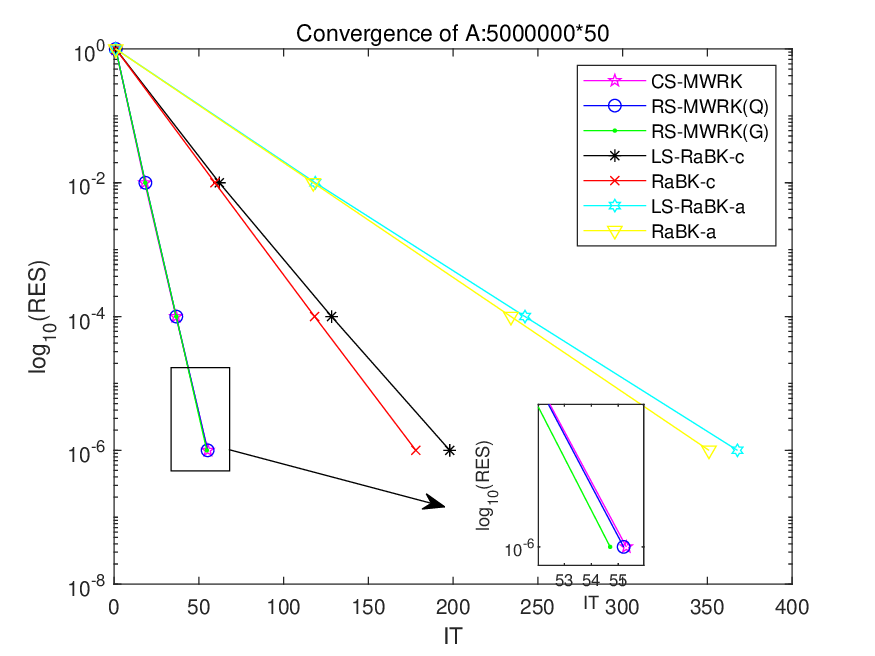}
	\includegraphics[width=0.4\linewidth,height=0.4\linewidth]{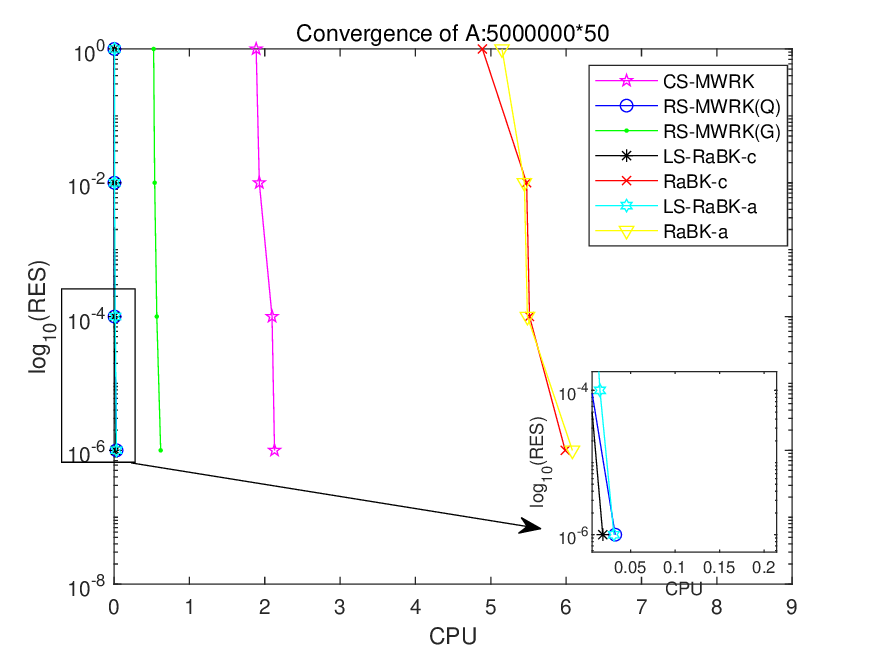}
	\caption{$log_{10}(RES)$ versus IT and CPU for CS-MWRK, RS-MWRK(Q), RS-MWRK(G), LS-RaBK(Q)-c, RaBK-c, LS-RaBK(Q)-a and RaBK-a when $A\in\mathbb{R}^{5000000\times 50}$.}
	\label{fig3}
\end{figure}
 Take  $A\in\mathbb{R}^{5000000\times 50}$, we draw Fig. \ref{fig3} to show the relationship between the error and the iterative steps(or the running time). Each line represents the median of the RES at that iteration or the CPU time after 50 trials. From the figures, we can see that  the RS-MWRK(G) method converges fastest under the same accuracy, while the LS-RaBK(Q)-c method requires the most petite runtime. It is noted that the running time does not start from $0$, the main reason is that the proposed algorithms require to compute ${\mathbf{S}}A$.
\clearpage
\section{Conclusions}\label{section6}
In this paper, we propose a new algorithmic framework based on the randomized matrices and two types of random matrices are applied into the framework. Besides, the convergence analysis are given. Furthermore, numerical experiments for the different sizes of the consistent linear systems are given to verify the effectiveness of the proposed algorithm, which show that the proposed algorithms are better than the original algorithms on a super large scale.
\section*{Acknowledgments}
The authors would like to thank the editor and reviewers for their
kind comments. This work was funded by National Natural Science Foundation of China (No. 12101136, 12201129), the Guangdong Basic and Applied Basic Research Foundations (Nos. 2023A1515011633), the Project of Science and Technology of Guangzhou (No. 2024A04J2056). All the authors thank for Center for Mathematics and Interdisciplinary Sciences, School of Mathematics and Statistics, Guangdong University of Technology.
\\
{\bf{Data Availability}} Data sharing not applicable to this article as no datasets were generated or analyzed during the current study.
\\
{\bf{Conflict of interest}} The authors declare that they have no conflict of interest.

\end{document}